\documentclass[11pt]{article}
\usepackage{amsfonts,amscd,amssymb, amsmath}
\newtheorem{definition}{Definition}[section]

\newtheorem{proposition}[definition]{Proposition}

\newtheorem{remark}[definition]{Remark}

\newtheorem{theorem}[definition]{Theorem}

\def\tl{\triangleleft}
\def\tr{\triangleright}

\newcommand{\nat}{\mbox{$\;\natural \;$}}

\def\rawo\lonra{\longrightarrow}

\def\ot{\otimes}

\allowdisplaybreaks[3]

\newenvironment{proof}{{\it Proof.}}{\hfill $ \square $ \vskip 4mm}

\begin{document}
\title{L-R-smash biproducts, double biproducts and a braided category of 
Yetter-Drinfeld-Long bimodules}
\author{Florin Panaite\thanks {Research
carried out while the first author was visiting the University of Antwerp, 
supported by a postdoctoral fellowship  
offered by FWO (Flemish Scientific Research Foundation). This author was 
also partially supported by the programme CEEX of the Romanian 
Ministry of Education and
Research, contract nr. 2-CEx06-11-20/2006.}\\
Institute of Mathematics of the 
Romanian Academy\\ 
PO-Box 1-764, RO-014700 Bucharest, Romania\\
e-mail: Florin.Panaite@imar.ro
\and 
Freddy Van Oystaeyen\\
Department of Mathematics and Computer Science\\
University of Antwerp, Middelheimlaan 1\\
B-2020 Antwerp, Belgium\\
e-mail: Francine.Schoeters@ua.ac.be}
\date{}
\maketitle

\begin{abstract}
Let $H$ be a bialgebra and $D$ an $H$-bimodule algebra and 
$H$-bicomodule coalgebra. We find sufficient conditions on $D$ for the  
L-R-smash product algebra and coalgebra structures on $D\otimes H$ to 
form a bialgebra (in this case we say that $(H, D)$ is an L-R-admissible 
pair), called L-R-smash biproduct. The Radford biproduct is a particular 
case, and so is, up to isomorphism, a double biproduct with trivial 
pairing. We construct a prebraided monoidal   
category ${\cal LR}(H)$, whose objects  
are $H$-bimodules $H$-bicomodules $M$ endowed with left-left and 
right-right Yetter-Drinfeld module as well as left-right and 
right-left Long module structures over $H$, with the property that, 
if $(H, D)$ is an L-R-admissible pair, then $D$ is a bialgebra in 
${\cal LR}(H)$.  
\end{abstract}
\section*{Introduction}
${\;\;\;\;}$
The L-R-smash product over a cocommutative Hopf algebra was introduced 
and studied in a series of papers  
\cite{b1}, \cite{b2}, \cite{b3}, \cite{b4}, with motivation and 
examples coming from the theory of deformation quantization. This 
construction was generalized in \cite{pvo} to the case of arbitrary 
bialgebras (even quasi-bialgebras), as follows:  
if $H$ is a bialgebra and $D$ is an  
$H$-bimodule algebra, 
the L-R-smash product $D\nat H$ is an associative  
algebra structure defined on $D\ot H$ by the 
multiplication rule  
\begin{eqnarray*}
&&(d\nat h)(d'\nat h')=(d\cdot h'_2)(h_1\cdot d')
\nat h_2h'_1, \;\;\;\forall \; d, d'\in D, \; 
h, h'\in H. 
\end{eqnarray*}
It was proved in \cite{pvo} that, if $H$ is moreover a Hopf algebra with 
bijective antipode, then $D\nat H$ is isomorphic to a diagonal 
crossed product $D\bowtie H$ as in \cite{bpvo}, \cite{hn}; this result was 
used in \cite{panvan} to give a very easy proof of the fact that two 
bialgebroids introduced independently in \cite{cm} and \cite{kadison} are 
actually isomorphic. 

The dual construction of the L-R-smash product was introduced also 
in \cite{pvo} under the name L-R-smash coproduct; this is a 
coassociative coalgebra $D\nat H$, where $D$ is an $H$-bicomodule coalgebra. 
A natural problem, not treated in \cite{pvo}, is to see under what 
conditions, for a given $H$-bimodule algebra $H$-bicomodule coalgebra $D$, 
the L-R-smash product and coproduct structures on $D\otimes H$ form a 
bialgebra. It seems to be difficult to obtain (nicely-looking) 
necessary and sufficient conditions on $D$ for this to happen. The aim of 
the present paper is to present a list of {\em sufficient conditions}, 
looking resonably nice and being general enough to cover some existing 
constructions from the literature. 

More precisely, if $D$ satisfies those conditions, we say that $(H, D)$ is 
an {\em L-R-admissible pair} and the bialgebra $D\nat H$ is called an 
{\em L-R-smash biproduct}. The Radford biproduct is a particular case, 
corresponding to the situation when the right action and coaction 
are trivial. We prove that a {\em double biproduct} $A\# H\# B$ 
(as in \cite{majid}, \cite{sommer}) with trivial pairing is isomorphic to 
an L-R-smash biproduct $(A\otimes B)\nat H$. Also, we show that a 
construction introduced in \cite{zhang} is a particular case of an 
L-R-smash biproduct. 

It is known that the Radford biproduct has a categorical 
interpretation (due to Majid): $(H, B)$ is an admissible pair (as in 
\cite{radf}) if and only if $B$ is a bialgebra in the Yetter-Drinfeld 
category $^H_H{\cal YD}$. We give a similar interpretation for  
L-R-admissible pairs. Namely, we define a prebraided category ${\cal LR}(H)$ 
(which is braided if $H$ has a skew antipode) consisting of $H$-bimodules 
$H$-bicomodules $M$ which are left-left and  
right-right Yetter-Drinfeld modules as well as left-right and 
right-left Long modules over $H$ (this category contains    
$^H_H{\cal YD}$ and ${\cal YD}_H^H$ as braided subcategories). We prove 
that all except one of the conditions for $(H, D)$ to be an L-R-admissible 
pair are equivalent to $D$ being a bialgebra in ${\cal LR}(H)$. The 
extra condition reads 
\begin{eqnarray*}
&&c^{<0>}\cdot d^{(-1)}\otimes c^{<1>}\cdot d^{(0)}=c\otimes d, 
\;\;\;\forall \;c, d\in D, 
\end{eqnarray*}        
and unfortunately does not seem to have a categorical interpretation 
inside ${\cal LR}(H)$. 
\section{The L-R-smash biproduct}
\setcounter{equation}{0}
${\;\;\;\;}$
We work over a field $k$. All algebras, linear spaces etc. will be 
over $k$; unadorned $\otimes $ means $\otimes _k$. For a bialgebra $H$ 
with comultiplication $\Delta $ we denote $\Delta (h)=h_1\otimes h_2$, 
for $h\in H$. For terminology concerning bialgebras, Hopf algebras and 
monoidal categories we refer to \cite{k}, \cite{m}.  

Let $H$ be a bialgebra and let $D$ be a vector space satisfying the 
following conditions: \\
(i) $D$ is an $H$-bimodule, with actions $h\otimes d\mapsto 
h\cdot d$ and $d\otimes h\mapsto d\cdot h$, for $h\in H$ and $d\in D$;\\
(ii) $D$ is an algebra, with unit $1_D$ and multiplication 
$c\otimes d\mapsto cd$, for $c, d\in D$; \\
(iii) $D$ is an $H$-bimodule algebra, that is $h\cdot 1_D=
\varepsilon (h)1_D$, $1_D\cdot h=\varepsilon (h)1_D$, 
$h\cdot (cd)=(h_1\cdot c)(h_2\cdot d)$ and $(cd)\cdot h=(c\cdot h_1)
(d\cdot h_2)$, for all $h\in H$ and $c, d\in D$; \\
(iv) $D$ is an $H$-bicomodule, with structures (for all $d\in D$):  
\begin{eqnarray*}
&&\rho :D\rightarrow D\ot H,\;\;\rho (d)=d^{<0>}\ot d^{<1>}, \\
&&\lambda :D\rightarrow H\ot D,\;\;\;\lambda (d)=d^{(-1)}\ot d^{(0)}; 
\end{eqnarray*}
(v) $D$ is a coalgebra, with comultiplication $\Delta _D:D\rightarrow 
D\ot D$, $\Delta _D(d)=d_1\ot d_2$, and counit $\varepsilon _D:D\rightarrow 
k$;\\
(vi) $D$ is an $H$-bicomodule coalgebra, that is, for all $d\in D$:
\begin{eqnarray*}
&&d_1^{(-1)}d_2^{(-1)}\ot d_1^{(0)}\ot d_2^{(0)}=d^{(-1)}\ot 
(d^{(0)})_1\ot (d^{(0)})_2, \label{lca} \\
&&d^{(-1)}\varepsilon _D(d^{(0)})=\varepsilon _D(d)1_H, \\
&&d_1^{<0>}\ot d_2^{<0>}\ot d_1^{<1>}d_2^{<1>}=(d^{<0>})_1\ot (d^{<0>})_2
\ot d^{<1>}, \label{rca} \\
&&\varepsilon _D(d^{<0>})d^{<1>}=\varepsilon _D(d)1_H.
\end{eqnarray*}

We denote the vector space $D\ot H$ by $D\nat H$ and elements 
$d\ot h$ by $d\nat h$. By \cite{pvo}, $D\nat H$ becomes an algebra (called 
L-R-smash product) with unit $1_D\nat 1_H$ and multiplication 
\begin{eqnarray*}
&&(d\nat h)(d'\nat h')=(d\cdot h'_2)(h_1\cdot d')\nat h_2h'_1, \;\;\;
\forall \;\;h, h'\in H, \;d, d'\in D, 
\end{eqnarray*}
and a coalgebra (called L-R-smash coproduct) with comultiplication 
and counit given by 
\begin{eqnarray*}
&&\Delta :D\nat H\rightarrow (D\nat H)\ot (D\nat H), \;\;\varepsilon :
D\nat H\rightarrow k, \\ 
&&\Delta (d\nat h)=(d_1^{<0>}\nat d_2^{(-1)}h_1)\ot (d_2^{(0)}\nat 
h_2d_1^{<1>}), \;\; 
\varepsilon (d\nat h)=\varepsilon _D(d)\varepsilon _H(h). 
\end{eqnarray*}

We consider now the following list of conditions, for $H$ and $D$ as above, 
corresponding to elements $h\in H$ and $c, d\in D$: 
\begin{eqnarray}
&&\varepsilon _D(1_D)=1, \;\;\;
\varepsilon _D(cd)=\varepsilon _D(c)\varepsilon _D(d), \label{eps1,2}\\
&&\varepsilon _D(h\cdot d)=\varepsilon _D(d\cdot h)=\varepsilon _D(d)
\varepsilon _H(h), \label{eps3}\\
&&\rho (1_D)=1_D\otimes 1_H, \;\;\;
\lambda (1_D)=1_H\otimes 1_D, \label{lrunit} \\
&&\Delta _D(1_D)=1_D\otimes 1_D, \label{deltaunit} \\
&&\rho (cd)=c^{<0>}d^{<0>}\otimes c^{<1>}d^{<1>}, \label{p1} \\
&&\lambda (cd)=c^{(-1)}d^{(-1)}\otimes c^{(0)}d^{(0)}, \label{p2} \\
&&\Delta _D(h\cdot d)=h_1\cdot d_1\otimes h_2\cdot d_2, \label{p3} \\
&&\Delta _D(d\cdot h)=d_1\cdot h_1\otimes d_2\cdot h_2, \label{p4} \\
&&\Delta _D(cd)=c_1(c_2^{(-1)}\cdot d_1^{<0>})\otimes 
(c_2^{(0)}\cdot d_1^{<1>})d_2, \label{co} \\
&&(h_1\cdot d)^{(-1)}h_2\otimes (h_1\cdot d)^{(0)}=
h_1d^{(-1)}\otimes h_2\cdot d^{(0)}, \label{y1} \\
&&(h\cdot d)^{<0>}\otimes (h\cdot d)^{<1>}=h\cdot d^{<0>}\otimes 
d^{<1>}, \label{y2} \\
&&(d\cdot h_2)^{<0>}\otimes h_1(d\cdot h_2)^{<1>}=
d^{<0>}\cdot h_1\otimes d^{<1>}h_2, \label{y3} \\
&&(d\cdot h)^{(-1)}\otimes (d\cdot h)^{(0)}=d^{(-1)}\otimes 
d^{(0)}\cdot h, \label{y4} \\
&&c^{<0>}\cdot d^{(-1)}\otimes c^{<1>}\cdot d^{(0)}=c\otimes d. \label{maj} 
\end{eqnarray}

If all these conditions hold, for all $h\in H$ and $c, d\in D$, by  
analogy with \cite{radf} we will say that $(H, D)$ is an 
{\bf L-R-admissible pair}.  
\begin{theorem} \label{main} 
If $(H, D)$ is an L-R-admissible pair, then $D\nat H$ with structures as 
above is a bialgebra,  
called the {\bf L-R-smash biproduct} of $D$ and $H$. 
\end{theorem}
\begin{proof}
It is very easy to see that $\varepsilon _{D\nat H}$ is an algebra map and 
$\Delta _{D\nat H}$ is unital, so we will only prove that 
$\Delta _{D\nat H}$ is multiplicative. We will prove first two 
auxiliary relations: 
\begin{eqnarray}
&&[c(h\cdot d)]_1\otimes [c(h\cdot d)]_2=c_1(c_2^{(-1)}h_1\cdot 
d_1^{<0>})\otimes (c_2^{(0)}\cdot d_1^{<1>})(h_2\cdot d_2), 
\label{cucu2} 
\end{eqnarray}
\begin{multline}
[c(h_1\cdot d)]_1\otimes [c(h_1\cdot d)]_2^{(-1)}h_2\otimes 
[c(h_1\cdot d)]_2^{(0)}=c_1(c_2^{(-1)}h_1\cdot d_1^{<0>})\otimes 
c_2^{(0)(-1)}h_2d_2^{(-1)}\\
\otimes (c_2^{(0)(0)}\cdot d_1^{<1>})(h_3\cdot d_2^{(0)}), \label{kaka}
\end{multline}  
for all $h\in H$, $c, d\in D$; we compute: 
\begin{eqnarray*}
[c(h\cdot d)]_1\otimes [c(h\cdot d)]_2&\overset{(\ref{co})}{=}&
c_1(c_2^{(-1)}\cdot (h\cdot d)_1^{<0>})\otimes 
(c_2^{(0)}\cdot (h\cdot d)_1^{<1>})(h\cdot d)_2\\
&\overset{(\ref{p3})}{=}&c_1(c_2^{(-1)}\cdot (h_1\cdot d_1)^{<0>})
\otimes (c_2^{(0)}\cdot (h_1\cdot d_1)^{<1>})(h_2\cdot d_2)\\
&\overset{(\ref{y2})}{=}&c_1(c_2^{(-1)}h_1\cdot d_1^{<0>})\otimes 
(c_2^{(0)}\cdot d_1^{<1>})(h_2\cdot d_2), 
\end{eqnarray*}
${\;\;\;\;\;}$$[c(h_1\cdot d)]_1\otimes [c(h_1\cdot d)]_2^{(-1)}h_2\otimes 
[c(h_1\cdot d)]_2^{(0)}$
\begin{eqnarray*}
&\overset{(\ref{cucu2})}{=}&c_1(c_2^{(-1)}h_1\cdot d_1^{<0>})\otimes 
[(c_2^{(0)}\cdot d_1^{<1>})(h_2\cdot d_2)]^{(-1)}h_3\otimes 
[(c_2^{(0)}\cdot d_1^{<1>})(h_2\cdot d_2)]^{(0)}\\
&\overset{(\ref{p2})}{=}&c_1(c_2^{(-1)}h_1\cdot d_1^{<0>})\otimes 
(c_2^{(0)}\cdot d_1^{<1>})^{(-1)}(h_2\cdot d_2)^{(-1)}h_3\otimes 
(c_2^{(0)}\cdot d_1^{<1>})^{(0)}(h_2\cdot d_2)^{(0)}\\
&\overset{(\ref{y1})}{=}&c_1(c_2^{(-1)}h_1\cdot d_1^{<0>})\otimes 
(c_2^{(0)}\cdot d_1^{<1>})^{(-1)}h_2d_2^{(-1)}\otimes  
(c_2^{(0)}\cdot d_1^{<1>})^{(0)}(h_3\cdot d_2^{(0)})\\
&\overset{(\ref{y4})}{=}&c_1(c_2^{(-1)}h_1\cdot d_1^{<0>})\otimes 
c_2^{(0)(-1)}h_2d_2^{(-1)}\otimes   
(c_2^{(0)(0)}\cdot d_1^{<1>})(h_3\cdot d_2^{(0)}), \;\;\;q.e.d.
\end{eqnarray*}
Let now $c, d\in D$ and $h, g\in H$; we compute:\\[2mm]
${\;\;\;}$
$\Delta ((c\nat h)(d\nat g))$
\begin{eqnarray*}
&=&\Delta ((c\cdot g_2)(h_1\cdot d)\nat h_2g_1)\\
&=&((c\cdot g_3)(h_1\cdot d))_1^{<0>}\otimes 
((c\cdot g_3)(h_1\cdot d))_2^{(-1)}h_2g_1\\
&&\otimes ((c\cdot g_3)(h_1\cdot d))_2^{(0)}\otimes 
h_3g_2((c\cdot g_3)(h_1\cdot d))_1^{<1>}\\
&\overset{(\ref{kaka})}{=}&[(c\cdot g_3)_1((c\cdot g_3)_2^{(-1)}h_1
\cdot d_1^{<0>})]^{<0>}\otimes (c\cdot g_3)_2^{(0)(-1)}h_2d_2^{(-1)}g_1\\
&&\otimes ((c\cdot g_3)_2^{(0)(0)}\cdot d_1^{<1>})(h_3\cdot d_2^{(0)})
\otimes h_4g_2[(c\cdot g_3)_1((c\cdot g_3)_2^{(-1)}h_1\cdot 
d_1^{<0>})]^{<1>}\\
&\overset{(\ref{p4}),\; (\ref{p1})}{=}& 
(c_1\cdot g_3)^{<0>}((c_2\cdot g_4)^{(-1)}h_1
\cdot d_1^{<0>})^{<0>}\otimes (c_2\cdot g_4)^{(0)(-1)}h_2d_2^{(-1)}g_1\\
&&\otimes ((c_2\cdot g_4)^{(0)(0)}\cdot d_1^{<1>})(h_3\cdot d_2^{(0)})
\otimes h_4g_2(c_1\cdot g_3)^{<1>}((c_2\cdot g_4)^{(-1)}h_1\cdot 
d_1^{<0>})^{<1>}\\
&\overset{(\ref{y2})}{=}&
(c_1\cdot g_3)^{<0>}((c_2\cdot g_4)^{(-1)}h_1
\cdot d_1^{<0><0>})\otimes (c_2\cdot g_4)^{(0)(-1)}h_2d_2^{(-1)}g_1\\
&&\otimes ((c_2\cdot g_4)^{(0)(0)}\cdot d_1^{<1>})(h_3\cdot d_2^{(0)})
\otimes h_4g_2(c_1\cdot g_3)^{<1>}d_1^{<0><1>}\\
&\overset{(\ref{y3})}{=}&
(c_1^{<0>}\cdot g_2)((c_2\cdot g_4)^{(-1)}h_1
\cdot d_1^{<0><0>})\otimes (c_2\cdot g_4)^{(0)(-1)}h_2d_2^{(-1)}g_1\\
&&\otimes ((c_2\cdot g_4)^{(0)(0)}\cdot d_1^{<1>})(h_3\cdot d_2^{(0)})
\otimes h_4c_1^{<1>}g_3d_1^{<0><1>}\\
&\overset{(\ref{y4})}{=}&
(c_1^{<0>}\cdot g_2)(c_2^{(-1)}h_1
\cdot d_1^{<0><0>})\otimes (c_2^{(0)}\cdot g_4)^{(-1)}h_2d_2^{(-1)}g_1\\
&&\otimes ((c_2^{(0)}\cdot g_4)^{(0)}\cdot d_1^{<1>})(h_3\cdot d_2^{(0)})
\otimes h_4c_1^{<1>}g_3d_1^{<0><1>}\\
&\overset{(\ref{y4})}{=}&
(c_1^{<0>}\cdot g_2)(c_2^{(-1)}h_1
\cdot d_1^{<0><0>})\otimes c_2^{(0)(-1)}h_2d_2^{(-1)}g_1\\
&&\otimes (c_2^{(0)(0)}\cdot g_4d_1^{<1>})(h_3\cdot d_2^{(0)})
\otimes h_4c_1^{<1>}g_3d_1^{<0><1>}, 
\end{eqnarray*}
\begin{eqnarray*}
\Delta (c\nat h)\Delta (d\nat g) 
&=&(c_1^{<0>}\nat c_2^{(-1)}h_1)(d_1^{<0>}\nat d_2^{(-1)}g_1)\otimes 
(c_2^{(0)}\nat h_2c_1^{<1>})(d_2^{(0)}\nat g_2d_1^{<1>})\\
&=&(c_1^{<0>}\cdot (d_2^{(-1)})_2g_2)((c_2^{(-1)})_1h_1\cdot d_1^{<0>})\nat 
(c_2^{(-1)})_2h_2(d_2^{(-1)})_1g_1\\
&&\otimes (c_2^{(0)}\cdot g_4(d_1^{<1>})_2)(h_3(c_1^{<1>})_1\cdot 
d_2^{(0)})\nat h_4(c_1^{<1>})_2g_3(d_1^{<1>})_1\\
&=&(c_1^{<0><0>}\cdot d_2^{(0)(-1)}g_2)(c_2^{(-1)}h_1\cdot 
d_1^{<0><0>})\otimes c_2^{(0)(-1)}h_2d_2^{(-1)}g_1\\
&&\otimes (c_2^{(0)(0)}\cdot g_4d_1^{<1>})(h_3c_1^{<0><1>}\cdot 
d_2^{(0)(0)})\otimes h_4c_1^{<1>}g_3d_1^{<0><1>}\\ 
&\overset{(\ref{maj})}{=}& 
(c_1^{<0>}\cdot g_2)(c_2^{(-1)}h_1
\cdot d_1^{<0><0>})\otimes c_2^{(0)(-1)}h_2d_2^{(-1)}g_1\\
&&\otimes (c_2^{(0)(0)}\cdot g_4d_1^{<1>})(h_3\cdot d_2^{(0)})
\otimes h_4c_1^{<1>}g_3d_1^{<0><1>}, 
\end{eqnarray*}
and we see that the two terms are equal.
\end{proof}
\begin{remark} {\em 
Obviously, the Radford biproduct (cf. \cite{radf}) is a particular case of 
the L-R-smash biproduct, corresponding to the case when the right action 
and coaction are trivial.} 
\end{remark}

We recall now from \cite{majid}, \cite{sommer} the construction of the 
so-called {\em double biproduct}, more precisely a particular case of it 
(corresponding to a trivial pairing, in the terminology of \cite{majid}). 
Let $H$ be a bialgebra, $A$ a bialgebra in the Yetter-Drinfeld category 
$_H^H{\cal YD}$ and $B$ a bialgebra in the Yetter-Drinfeld category 
${\cal YD}_H^H$, with the following notation for the structure maps: 
counits $\varepsilon _A$ and $\varepsilon _B$, comultiplications 
$\Delta _A(a)=a_1\otimes a_2$ and $\Delta _B(b)=b_1\otimes b_2$, and 
actions and coactions 
\begin{eqnarray*}
&&H\otimes A\rightarrow A, \;\;\;h\otimes a\mapsto h\triangleright a, \\
&&A\rightarrow H\otimes A, \;\;\;a\mapsto a^1\otimes a^2, \\
&&B\otimes H\rightarrow B, \;\;\;b\otimes h\mapsto b\triangleleft h, \\
&&B\rightarrow B\otimes H, \;\;\;b\mapsto b^1\otimes b^2, 
\end{eqnarray*}
for all $h\in H$, $a\in A$, $b\in B$. We denote by $A\# H\# B$ the vector 
space $A\otimes H\otimes B$ (the element $a\otimes h\otimes b$ is denoted 
by $a\# h\# b$),  
which becomes an algebra (called two-sided smash product) with unit 
$1_A\# 1_H\# 1_B$ and multiplication 
\begin{eqnarray*}
&&(a\# h\# b)(a'\# h'\# b')=a(h_1\tr a')\# h_2h'_1\# (b\tl h'_2)b', 
\end{eqnarray*} 
and a coalgebra (called two-sided smash coproduct) with counit 
$\varepsilon (a\# h\# b)=\varepsilon _A(a)\varepsilon _H(h)
\varepsilon _B(b)$ and comultiplication 
\begin{eqnarray*}
&&\Delta :A\# H\# B\rightarrow (A\# H\# B)\otimes (A\# H\# B), \\
&&\Delta (a\# h\# b)=(a_1\# a_2^1h_1\# b_1^1)\otimes 
(a_2^2\# h_2b_1^2\# b_2). 
\end{eqnarray*}
\begin{proposition} (\cite{majid}, \cite{sommer})  
Assume that moreover the following condition holds:
\begin{eqnarray}
&&b^2\tr a^2\otimes b^1\tl a^1=a\otimes b, \;\;\;\forall \;\;
a\in A, \;b\in B.  \label{som}
\end{eqnarray}
Then $A\# H\# B$ is a bialgebra, called the {\bf double biproduct}. 
\end{proposition}
\begin{proposition} \label{relat}
Let $A\# H\# B$ be a double biproduct bialgebra. Define $D=A\otimes B$, 
with tensor product algebra and coalgebra structures and with two-sided 
actions and coactions given by  
\begin{eqnarray*}
&&H\otimes (A\otimes B)\rightarrow A\otimes B, 
\;\;\;h\otimes (a\otimes b)\mapsto h\cdot (a\otimes b):= 
h\triangleright a\otimes b, \\
&&A\otimes B\rightarrow H\otimes (A\otimes B), \;\;\;a\otimes b
\mapsto (a\otimes b)^{(-1)}\otimes (a\otimes b)^{(0)}:= 
a^1\otimes (a^2\otimes b), \\
&&(A\otimes B)\otimes H\rightarrow A\otimes B, \;\;\;(a\otimes b)\otimes h
\mapsto (a\otimes b)\cdot h:= 
a\otimes b\triangleleft h, \\
&&A\otimes B\rightarrow (A\otimes B)\otimes H, \;\;\;a\otimes b\mapsto 
(a\otimes b)^{<0>}\otimes (a\otimes b)^{<1>}:=
(a\otimes b^1)\otimes b^2. 
\end{eqnarray*}
Then $(H, D)$ is an L-R-admissible pair and we have a  
bialgebra isomorphism 
\begin{eqnarray*}
&&\phi :(A\otimes B)\nat H\simeq A\# H\# B, \;\;\;(a\otimes b)\nat h
\mapsto a\# h\# b. 
\end{eqnarray*} 
\end{proposition}
\begin{proof}
The fact that $(H, D)$ is an L-R-admissible pair follows  
by direct computation; let us only check (\ref{maj}), for $a, a'\in A$ and 
$b, b'\in B$: 
\begin{eqnarray*}
(a\otimes b)^{<0>}\cdot (a'\otimes b')^{(-1)}\otimes (a\otimes b)^{<1>}
\cdot (a'\otimes b')^{(0)}&=&(a\otimes b^1)\cdot a'^1\otimes 
b^2\cdot (a'^2\otimes b')\\
&=&(a\otimes b^1\tl a'^1)\otimes (b^2\tr a'^2\otimes b')\\
&\overset{(\ref{som})}{=}&(a\otimes b)\otimes (a'\otimes b'), \;\;\;q.e.d. 
\end{eqnarray*}
We know from \cite{pvo}, Proposition 2.4, that $\phi $ is an algebra 
isomorphism, and an easy computation shows that $\phi $ is also a 
coalgebra map. 
\end{proof}

We recall now the following result from \cite{zhang}. Let $H$ be a 
bialgebra and $D$ an $H$-bimodule bialgebra (i.e. $D$ is a bialgebra which 
is an $H$-bimodule algebra and an $H$-bimodule coalgebra). Consider the 
L-R-smash product algebra $D\nat H$, together with the tensor product 
coalgebra structure on it (i.e. $\Delta (d\nat h)=(d_1\nat h_1)\otimes 
(d_2\nat h_2)$ and $\varepsilon (d\nat h)=\varepsilon _D(d)
\varepsilon _H(h)$). Then $D\nat H$ with these structures is a bialgebra 
if and only if the following conditions are satisfied, for all 
$h\in H$, $d\in D$: 
\begin{eqnarray}
&&h_1\cdot d\otimes h_2=h_2\cdot d\otimes h_1, \label{LRS1}\\
&&d\cdot h_1\otimes h_2=d\cdot h_2\otimes h_1. \label{LRS2}
\end{eqnarray}
    
This result is a particular case of Theorem \ref{main}. Indeed, consider on 
$D$ the left and right trivial $H$-coactions (i.e. $d^{(-1)}\otimes 
d^{(0)}=1_H\otimes d$ and $d^{<0>}\otimes d^{<1>}=d\otimes 1_H$, 
for $d\in D$). Then one can easily check that $(H, D)$ is an 
L-R-admissible pair ((\ref{LRS1}) and (\ref{LRS2}) are precisely (\ref{y1})  
and respectively (\ref{y3})) and the L-R-smash coproduct coalgebra 
structure in this case coincides with the tensor product coalgebra 
structure.   
\section{A braided category related to L-R-smash biproducts}
\setcounter{equation}{0}
${\;\;\;\;}$Let $H$ be a bialgebra. We will introduce a prebraided category 
associated to $H$, denoted by ${\cal LR}(H)$. The objects of 
${\cal LR}(H)$ are vector spaces $M$ endowed with $H$-bimodule and 
$H$-bicomodule structures (denoted by $h\otimes m\mapsto h\cdot m$, 
$m\otimes h\mapsto m\cdot h$, $m\mapsto m^{(-1)}\otimes m^{(0)}$, 
$m\mapsto m^{<0>}\otimes m^{<1>}$, for all $h\in H$, $m\in M$), such that 
$M$ is a left-left Yetter-Drinfeld module, a left-right Long module, 
a right-right Yetter-Drinfeld module and a right-left Long module, i.e.  
\begin{eqnarray}
&&(h_1\cdot m)^{(-1)}h_2\otimes (h_1\cdot m)^{(0)}=
h_1m^{(-1)}\otimes h_2\cdot m^{(0)}, \label{caty1} \\
&&(h\cdot m)^{<0>}\otimes (h\cdot m)^{<1>}=h\cdot m^{<0>}\otimes 
m^{<1>}, \label{caty2} \\
&&(m\cdot h_2)^{<0>}\otimes h_1(m\cdot h_2)^{<1>}=
m^{<0>}\cdot h_1\otimes m^{<1>}h_2, \label{caty3} \\
&&(m\cdot h)^{(-1)}\otimes (m\cdot h)^{(0)}=m^{(-1)}\otimes 
m^{(0)}\cdot h, \label{caty4}
\end{eqnarray}
for all $h\in H$, $m\in M$. The morphisms in ${\cal LR}(H)$ are the 
$H$-bilinear $H$-bicolinear maps. 

One can check that ${\cal LR}(H)$ becomes a strict monoidal category, 
with unit $k$ endowed with usual $H$-bimodule and $H$-bicomodule structures, 
and tensor product given as follows: if $M, N\in {\cal LR}(H)$ then 
$M\otimes N\in {\cal LR}(H)$ with structures (for all $m\in M$, $n\in N$, 
$h\in H$): 
\begin{eqnarray*}
&&h\cdot (m\otimes n)=h_1\cdot m\otimes h_2\cdot n, \\
&&(m\otimes n)\cdot h=m\cdot h_1\otimes n\cdot h_2, \\
&&(m\otimes n)^{(-1)}\otimes (m\otimes n)^{(0)}=m^{(-1)}n^{(-1)}\otimes 
(m^{(0)}\otimes n^{(0)}), \\
&&(m\otimes n)^{<0>}\otimes (m\otimes n)^{<1>}=(m^{<0>}\otimes n^{<0>})
\otimes m^{<1>}n^{<1>}.
\end{eqnarray*}
\begin{proposition} The monoidal category 
${\cal LR}(H)$ is a prebraided category, with braiding defined, for all
$M, N\in {\cal LR}(H)$, $m\in M$, $n\in N$, by 
\begin{eqnarray*}
&&c_{M, N}:M\otimes N\rightarrow N\otimes M, \;\;\;
c_{M, N}(m\otimes n)=m^{(-1)}\cdot n^{<0>}\otimes m^{(0)}\cdot n^{<1>}.
\end{eqnarray*}
If $H$ has a skew antipode $S^{-1}$, then ${\cal LR}(H)$ is 
braided, the inverse of $c$ being given by 
\begin{eqnarray*}
&&c_{M, N}^{-1}:N\otimes M\rightarrow M\otimes N, \;\;\;
c_{M, N}^{-1}(n\otimes m)=m^{(0)}\cdot S^{-1}(n^{<1>})\otimes 
S^{-1}(m^{(-1)})\cdot n^{<0>}. 
\end{eqnarray*}
\end{proposition}
\begin{proof}
We only check that $c$ is left $H$-linear, right $H$-colinear and 
satisfies one of the two hexagonal equations, and leave the rest to the 
reader. For $M, N, P\in {\cal LR}(H)$ and $h\in H$, $m\in M$, $n\in N$, 
$p\in P$, we compute: 
\begin{eqnarray*}
c_{M, N}(h\cdot (m\otimes n))&=&c_{M, N}(h_1\cdot m\otimes h_2\cdot n)\\
&=&(h_1\cdot m)^{(-1)}\cdot (h_2\cdot n)^{<0>}\otimes (h_1\cdot m)^{(0)}
\cdot (h_2\cdot n)^{<1>}\\
&\overset{(\ref{caty2})}{=}&
(h_1\cdot m)^{(-1)}h_2\cdot n^{<0>}\otimes (h_1\cdot m)^{(0)}
\cdot n^{<1>}\\
&\overset{(\ref{caty1})}{=}&
h_1m^{(-1)}\cdot n^{<0>}\otimes h_2\cdot m^{(0)}
\cdot n^{<1>}\\
&=&h_1\cdot (m^{(-1)}\cdot n^{<0>})\otimes h_2\cdot (m^{(0)}
\cdot n^{<1>})\\
&=&h\cdot c_{M, N}(m\otimes n), 
\end{eqnarray*} 
${\;\;\;\;\;}$
$(\rho _{N\otimes M}\circ c_{M, N})(m\otimes n)$
\begin{eqnarray*}
&=&\rho _{N\otimes M}(m^{(-1)}\cdot n^{<0>}\otimes m^{(0)}\cdot n^{<1>})\\
&=&(m^{(-1)}\cdot n^{<0>})^{<0>}\otimes (m^{(0)}\cdot n^{<1>})^{<0>}\otimes 
(m^{(-1)}\cdot n^{<0>})^{<1>}(m^{(0)}\cdot n^{<1>})^{<1>}\\
&\overset{(\ref{caty2})}{=}&
m^{(-1)}\cdot n^{<0><0>}\otimes (m^{(0)}\cdot n^{<1>})^{<0>}\otimes 
n^{<0><1>}(m^{(0)}\cdot n^{<1>})^{<1>}\\
&=&m^{(-1)}\cdot n^{<0>}\otimes (m^{(0)}\cdot (n^{<1>})_2)^{<0>}\otimes 
(n^{<1>})_1(m^{(0)}\cdot (n^{<1>})_2)^{<1>}\\
&\overset{(\ref{caty3})}{=}&
m^{(-1)}\cdot n^{<0>}\otimes m^{(0)<0>}\cdot (n^{<1>})_1\otimes  
m^{(0)<1>}(n^{<1>})_2\\
&=&m^{<0>(-1)}\cdot n^{<0><0>}\otimes m^{<0>(0)}\cdot n^{<0><1>}\otimes   
m^{<1>}n^{<1>}\\
&=&c_{M, N}(m^{<0>}\otimes n^{<0>})\otimes m^{<1>}n^{<1>}\\
&=&(c_{M, N}\otimes id_H)\circ \rho _{M\otimes N}(m\otimes n), 
\end{eqnarray*}
${\;\;\;\;\;}$
$(id_N\otimes c_{M, P})\circ (c_{M, N}\otimes id_P)(m\otimes n\otimes p)$
\begin{eqnarray*}
&=&(id_N\otimes c_{M, P})(m^{(-1)}\cdot n^{<0>}\otimes m^{(0)}\cdot 
n^{<1>}\otimes p)\\
&=&m^{(-1)}\cdot n^{<0>}\otimes (m^{(0)}\cdot n^{<1>})^{(-1)}\cdot 
p^{<0>}\otimes (m^{(0)}\cdot n^{<1>})^{(0)}\cdot p^{<1>}\\
&\overset{(\ref{caty4})}{=}&
m^{(-1)}\cdot n^{<0>}\otimes m^{(0)(-1)}\cdot  
p^{<0>}\otimes m^{(0)(0)}\cdot n^{<1>}p^{<1>}\\
&=&(m^{(-1)})_1\cdot n^{<0>}\otimes (m^{(-1)})_2\cdot   
p^{<0>}\otimes m^{(0)}\cdot n^{<1>}p^{<1>}\\
&=&m^{(-1)}\cdot (n\otimes p)^{<0>}\otimes m^{(0)}\cdot (n\otimes p)^{<1>}\\
&=&c_{M, N\otimes P}(m\otimes n\otimes p).
\end{eqnarray*}
Also, the bijectivity of $c$ in the presence of a skew antipode follows 
by a direct computation which is left to the reader. 
\end{proof}
\begin{remark}{\em 
We denote as usual by $^H_H{\cal YD}$ and ${\cal YD}_H^H$ the categories 
of left-left and respectively right-right Yetter-Drinfeld modules over $H$. 
One can check that, if $V\in $$\;^H_H{\cal YD}$ and $W\in {\cal YD}_H^H$, 
then $V\otimes W\in {\cal LR}(H)$, with structures as in Proposition 
\ref{relat}. In particular, for $W=k$ and respectively $V=k$, we obtain 
that $^H_H{\cal YD}$ and ${\cal YD}_H^H$ are subcategories of 
${\cal LR}(H)$, and one can see that they are actually braided subcategories, 
i.e. the braiding of ${\cal LR}(H)$ restricts to the usual braidings of 
$^H_H{\cal YD}$ and ${\cal YD}_H^H$. }
\end{remark}

We can state now the categorical interpretation of L-R-admissible pairs: 
\begin{proposition}
Let $H$ be a bialgebra and $D$ a vector space. Then $(H, D)$ is an 
L-R-admissible pair if and only if $D$ is a bialgebra in ${\cal LR}(H)$ 
satisfying (\ref{maj}).  
\end{proposition}
\begin{proof}
A straightforward verification; we only note that (\ref{co}) expresses 
the fact that the comultiplication of $D$ is an algebra map inside the 
category ${\cal LR}(H)$. 
\end{proof}

\end{document}